# On planar regular graphs degree three without Hamiltonian cycles[1]

E. Grinbergs

Computing Centre of Latvian State University

**Abstract.** Necessary condition to have Hamiltonian cycle in planar graph is given. Examples of regular planar graphs degree three without Hamiltonian cycle are built.

In this article, by 'graph' we mean finite connected nonoriented graph without loops, multiedges, and vertices of degree one.

Let us call *the nontrivial cut of the graph* any set of edges removing which graph falls into two connected components with at least three vertices in each. The minimal number of edges contained in the nontrivial cut we call the *quasiconnectivity*. Finally, finite regular planar graph degree three we call the *map* if its quasiconnectivity is at least two, i.e. it does not have isthmi. It is well known that map may not have Hamiltonian cycle. However, for the only example of the map without Hamiltonian cycle that repeatedly appears in literature [1, p. 210; 2, p. 57], quasiconnectivity is but three (with besides three minimal nontrivial cuts).

On the other side, mostly interesting are just these maps that have quasi-connectivity of its maximal value -- five, at least as the only possible counterexamples of the hypothesis of the four colour theorem. We are going to establish one characteristic of planar graphs having Hamiltonian cycle that would, in its turn, let us build maps of quasiconnectivity five or four without such cycles.

Let us consider planar embedding of the graph $G$ (i.e. planar topological graph according Berge, depicting $G$) with definite Hamiltonian cycle $H$ in it. The edges not contained in $H$ let us call the *chords*. Because of the planarity of considered embedding, $H$ is a closed Jordan curve without selfcrossings and, consequently, divides the plane into two parts -- outer, containing infinitely remote point of the plane, -- and inner. Each chord has end vertices on $H$ passing only along outer or inner part of the plane. These chords, together with faces into which the two parts of the plane are divided, may be called outer and inner respectively.

We define the *weight* of the face bordered by elementary cycle to be the length of the bordering cycle [or number of vertices of the border] minus two. If $\nu$ faces are considered so that the sum of the length of bordering cycles of them is equal to $\Sigma$, the sum $S$ of their weights may be expressed as follows:
$$S = \Sigma - 2\nu. \qquad (1)$$

Let now $h$ -- length of the considered Hamiltonian cycle H. $S_1$ and $S_2$ are sums of weights of inner and outer faces respectively. Then
$$S_1 = S_2 = h - 2. \qquad (2)$$

---

To prove it, let us eliminate all chords. Then equalities hold in (2) by definition. Let us now restore all chords one by one. By restoring each subsequent chord, addends in $S_1$ and $S_2$, and sums itself, do not change. For the first sum, let us follow changes in the right side of (1). Number $v$ increases by one because one of the faces, say $A$, due to added chord, partitions into two new faces. Sum $\Sigma$ increases by two because all edges constituting border of $A$ go also into border of just one of the new faces. Thus, value of the left side does not change. In consequence, expression (2) remains true after all chords are restored in their place.

If the set of all faces of an embedding of the graph in plane are partitioned into two subsets with equal total weights $S_1$ and $S_2$, i.e.
$$S_1 = S_2, \tag{3}$$
then such partition may be called *isobaric* but the set of all edges contained in both these partitions -- the *border* of that partition. Then the obtained result may be formulated by the following theorem:

Theorem. *If planar graph $G$ has Hamiltonian cycle $H$, then for any embedding in plane there exists an isobaric partition so that the border of it is $H$.*

Let us add some remarks. Firstly, it is clear that the formulated theorem remains true for planar multigraphs too (of course, without loops, isthmi, and vertices of first order). Secondly, for our aims important is just expression (3) because expression
$$S_1 + S_2 = 2(h-2)$$
following from second sign of equality in (2) is simply some alternative form of Euler's formula. Finally, non-existence of isobaric partition may be rather convenient sufficient condition for non-existence of Hamiltonian cycle. However, if such partitions are at hand, it is hardly purposeful to use these conditions of the formulated theorem to solve the problem about existence of the Hamiltonian cycle in the given planar graph or to find such cycles.

The theorem should be useful the other way around: if at first, the set of available weights of edges is limited, then it is possible to build planar embeddings that one of the following holds:
a) isobaric partition does not exist;
b) isobaric partition exists, but no one of them has Hamiltonian cycle as its boundary.

Acquired planar graphs, by virtue of the theorem, should be without Hamiltonian cycles. Maps satisfying a) may be built if for its building, it is used just one face with the weight not multiple of three, and arbitrary number of faces of weights that are multiples of three, i.e. faces with $3k+2$ vertices for natural $k$. Then, for any partition of the set of all faces into two subsets, only for one part sum of weights should be multiple of three. Isobaric partitions, as follows, should not be there.

In order to build map with b) holding, we may take three faces with common vertex with weights that are comparable by module of three, but not multiples of three, but all other faces, as higher, with weights multiples of three. Then, for any isobaric partition, first three faces should go into one subset of partition; border of partition would not go through $x$, and, as follows, Hamiltonian cycle would not be there.

It seems to us that the building of the map $G$ is easer to perform for the dual graph of the map we are looking for, i.e. for planar triangulation of $G'$. For Hamiltonian cycle $H$ in $G$, there corresponds the cut $H'$ in $G'$ that separates one from another two trees, which may be taken as bichromatic components to get proper colouring of vertices of $G'$. However, if in $G$ there is no Hamiltonian cycle, then $G'$ is without such cuts; for proper colouring of vertices of $G'$, any bichromatic component is not tree but either has even cycle or isn't connected or both. For the length of boundary cycles of faces of $G$, there correspond degrees of vertices in $G'$ so that the task can be reduced to finding planar triangulations with arbitrary number of vertices of degree 5, 8, 11, ... and one vertex of degree not comparable with two by module three in case a), or three pairwise adjacent vertices of degree comparable between themselves but not two by module three in case b).

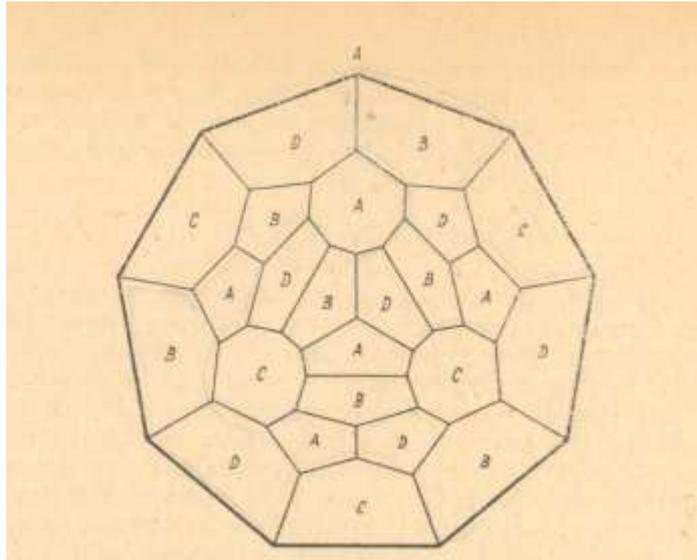

Fig. 1. Map without Hamiltonian cycles:
$f_9 = 1; f_5 = 21; f_8 = 3; f = 25; q = 5$

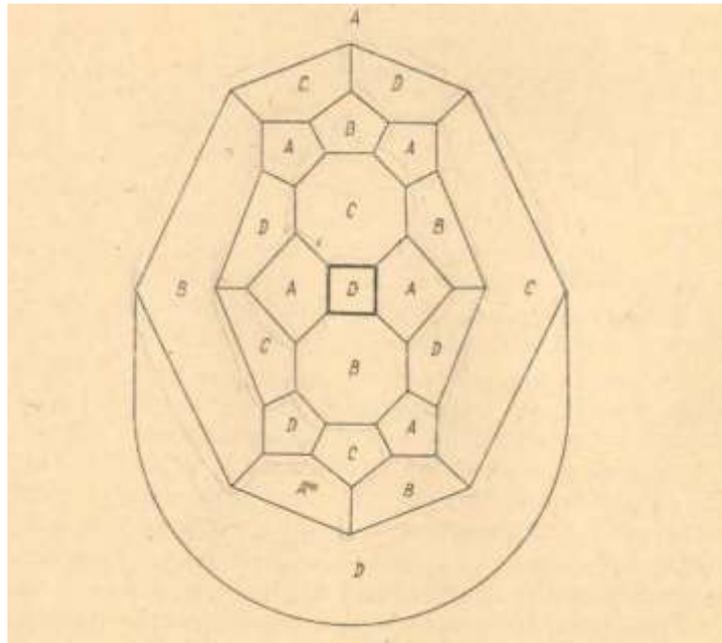

Fig. 2. Map without Hamiltonian cycles:
$f_4 = 1; f_5 = 18; f_8 = 4; f = 23; q = 4$

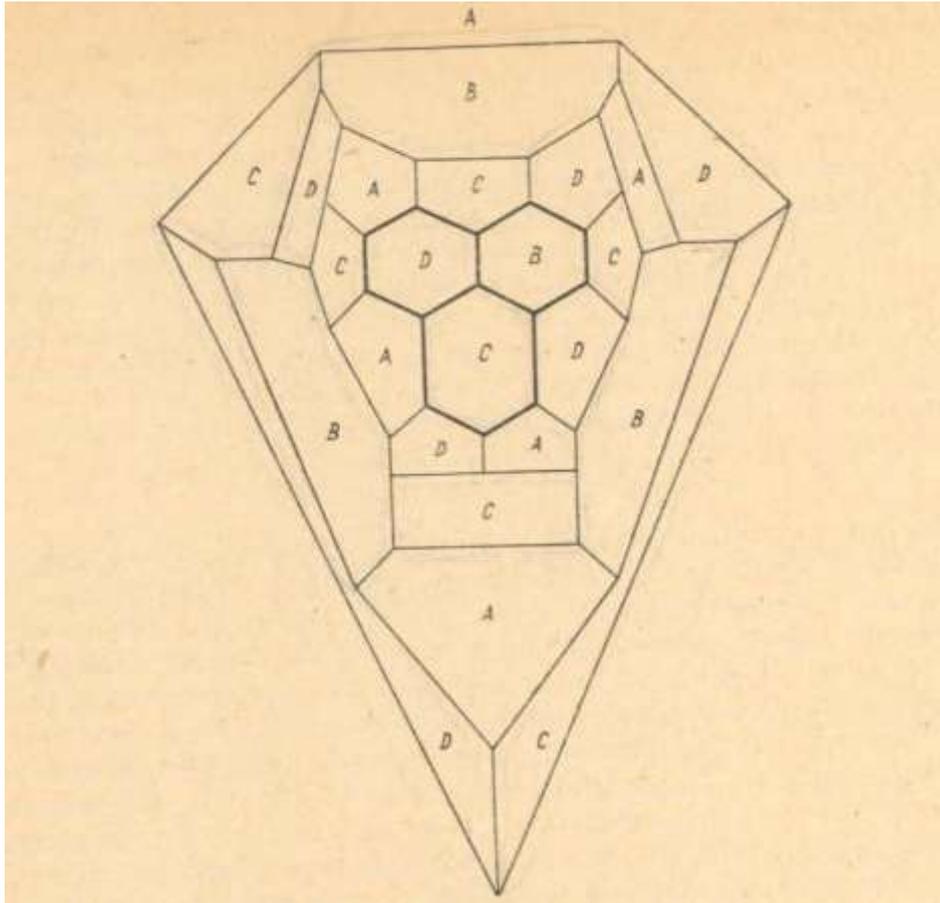

Fig. 3. Map without Hamiltonian cycles:
$f_6 = 3; f_5 = 18; f_8 = 3; f = 24; q = 5$

On fig. 1-3, some maps obtained in the described way are given. Borders of 'special' faces, weights of which are not multiples of three, are drawn with bold lines. We will denote by $f_i$ -- number of faces with $i$ vertices, by $f$ - summary number of faces; by $q$ -- quasi-connectivity; by $A, B, C, D$ -- colours of one of the proper colourings of faces with four colours. The map given on fig. 1 is dual to triangulation on fig. 1, I in [3]. Maps shown on fig. 1 and 2 correspond to case a) but on 3 -- case b). Among the maps built by us in the described way, two last maps (fig. 2, 3) have least possible $f$ for $q = 4$ and $5$. It would be interesting to state minimal values for maps without Hamiltonian cycles with $q = 3, 4$ or $5$. Since maps from [1,2] have $f = 25$ at $q = 3$, upper limits for these minimal values are known.

For the general number of faces, edges, vertices, or maximal number of vertices of one face, maps without Hamiltonian cycle does not have upper limit. Really, it is possible to show the way of building maps with each of these values being arbitrary large. Let us describe one way of such building for dual triangulation $G'$.

We choose natural numbers $\alpha$ and $\beta$ with condition
$$\beta = 2 \pmod 3 \tag{4}$$
holding. We are building elementary cycle $C_1$ of length $3^\alpha \beta$ with another elementary cycle $C_0$ within it of the same length, within which we put vertex $x$ connecting it via edges with all vertices of $C_0$. After that, in a zigzag way, we triangulize area between $C_0$ and $C_1$: we connect by edge one vertex of $C_0$ with one vertex of $C_1$, then after, we build new edges by moving end points alternatively along $C_0$ and $C_1$ by one edge in the same direction. In result, we get triangulation of the inside of $C_1$; vertex $x$ has degree $3^\alpha \beta$, all vertices of $C_0$ have degree 5, but to vertices of $C_1$ from inside, there come edges in twos. Further, by $\alpha$ steps of the same type we build parts of $G'$ that lay between cycles $C_i$ and $C_{i+1}$ ($i = 1, 2, ..., \alpha$), namely, we partition all vertices of $C_i$ into triples of following one after another vertices. To vertices of each triple we assign degrees 8, 5, 8, and in accordance with layout of fig. 4, we build transition to $C_{i+1}$ (the numbers indicated are degrees of vertices).

By feasibility of splitting vertices in triples and by getting prescribed values of degrees, following features of our cycles are ensured:
1) cycles of $C_i$ have lengths $3^{\alpha-i+1} 4^{i-1} \beta$;

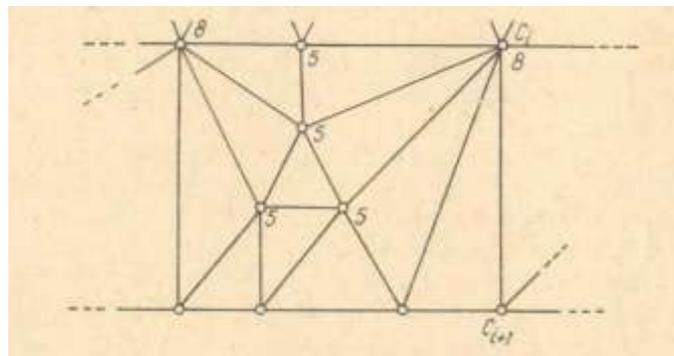

Fig. 4 Layout for building graph $G'$.

2) from inside to vertices of $C_i$ come two edges. These features are true because they hold firstly for $i=1$ and further, when we move from $i$ to $i+1$, because the construction being used retains them during augmentation.

We end the building of the triangulation by uniting all vertices of $C_{\alpha+1}$ with vertex $z$ taken outside $C_{\alpha+1}$. The degree of $z$ should be equal to $4^{\alpha}\beta$, and by virtue of (4), we have:
$$4^{\alpha}\beta \equiv \beta \equiv 2 (\text{mod } 3) \ .$$

Now, if we go over to map $G_i$ that is dual to $G'$, then face corresponding to $x$ should have $3^{\alpha}\beta$ vertices, consequently, weight $3^{\alpha}\beta - 2$ that is not multiple of three. The weights of all other faces would be 3, 6 or $4^{\alpha}\beta - 2$, and therefore multiples of three so that no isobaric partitions, no Hamiltonian cycles are possible.

It is easy to see that all graphs $G'$ built in the described way are four-chromatic.

Finally, let us remark that by stated theorem it is easy to get analogue proposition for maps that are 3H-graphs. As is well known, regular graph degree three is called 3H-graph if there exists such regular colouring of edges with three colours that edges of each two colours compose Hamiltonian cycle.

Let map $G$ be 3H-graph and $H_1, H_2, H_3$ be Hamiltonian cycles by virtue of definition of 3H-graph. The following are known facts that may be easily proven:

- we get regular colouring of faces of $G$ with colours 0,1,2,3 if number of colour of each face equals number of cycles $H_i$ within which this face lies;

- by suitable numbering of colours of edges present in the definition of 3H-graph, regular colouring of edges with colours 0, 1, 2 is obtained if edge that touches faces with colours $i$ and $j$ obtains colour numbered $|i+j-3|$. Let, finally, $\sigma_i$ be sum of weights of all faces with colour $i$. Using the theorem for each Hamiltonian cycle $H_j$, we see that sum of each two values of $\sigma_i$ equals to sum of both others, consequently, all $\sigma_i$ have the same value.

Thus, the following holds:

Corollary. If map is 3H-graph, then exists such partition of the set of all faces into four subsets $A_i (i=0,1,2,3)$ having equal summary weights that colouring of all faces of subset $A_i$ with $i$ colours is regular and have all indicated above features.

Exempli gratia of a map that is 3H-graph, it is possible to point to the projection of dodecahedron that has all faces of weight three.

E. Grinberg

**On Planar Gomogenous Graphs Degree Three Without Hamilton Circuits[2]**

**Annotation**

A necessary condition for Hamilton circuits of a planar graph is given in this paper. There are given examples of gomogenous planar graphs degree three having no Hamilton circuits.

Submitted on 2 August 1967.

---

[2] This annotation is given in the form it was published in 1968.